\documentclass[12pt]{amsart}
\newcommand{\Q}{\mathbb{Q}}
\newcommand{\F}{\mathbb{F}}

\newcommand{\Z}{\mathbb{Z}}
\renewcommand{\P}{\mathbb{P}}
\newcommand{\rank}{{\rm rank}~}
\usepackage{amsthm,url}
\newtheorem{thm}{Theorem}
\newtheorem{lem}[thm]{Lemma}
\newtheorem*{rem}{Remark}

\newtheorem*{ack}{Acknowledgements}
\usepackage{amssymb}
\usepackage{comment}
\usepackage{hyperref}
\begin{document}

\title{Rational points on $x^{3} + x^{2} y^{2} + y^{3} = k$}
\author{Xiaoan Lang}
\author{Jeremy Rouse}
\begin{abstract}
  We study the problem of determining, given an integer $k$, the rational solutions to $C_{k} : x^{3}z + x^{2} y^{2} + y^{3}z = kz^{4}$. For $k \ne 0$, the curve $C_{k}$ has genus $3$ and there are maps from $C_{k}$ to three elliptic curves $E_{1,k}$, $E_{2,k}$, $E_{3,k}$. We explicitly determine the rational points on $C_{k}$ under the assumption that one of these elliptic curves has rank zero. We discuss the challenges involved in extending our result to handle all $k \in \mathbb{Q}$.
\end{abstract}

\maketitle

\section{Introduction and Statement of Results}
\label{intro}

A fundamental problem in number theory is to determine the rational solutions to an equation or system of equations with rational coefficients. In 1983,
Gerd Faltings proved that if $C$ is a curve of genus at least $2$ defined
over a number field $K$, then $C(K)$ is finite. However, the proof is not constructive, and it is still not known if there is an algorithm that will determine
the finite set $C(K)$ for any curve $C$ of genus at least $2$ and any number field $K$. 

A variety of techniques exist that in some cases allow one to determine the set $C(K)$. The most basic is that if $C$ does not have local points at some completion of $K$, then $C(K) = \emptyset$. The Chabauty-Coleman method provides
a way to determine $C(K)$ under the assumption that the Mordell-Weil rank of the Jacobian is less than the genus. There are also a handful of other techniques that work in more restricted situations (including maps to lower genus curves, \'etale descent, and elliptic curve Chabauty). In addition, there is the relatively new Chabauty-Kim method.

In this paper, we consider the family of genus $3$ curves $C_{k} : x^{3}z + x^{2} y^{2} + y^{3}z = kz^{4}$ over $\Q$ under the restriction that $k \in \mathbb{Z}$.
For each $k$, this curve has two obvious rational points: $(1 : 0 : 0)$ and $(0 : 1 : 0)$. The map which interchanges $x$ and $y$ is an automorphism of $C_{k}$ over $\Q$, and generically the automorphism group of $C_{k}/\overline{\Q}$ is isomorphic to $S_{3}$. As a consequence, the curve admits three independent maps to elliptic curves. In particular, define
\begin{align*}
  E_{1,k} &: y^{2} + 3xy = x^{3} + k\\
  E_{2,k} &: y^{2} = x^{3} + 4kx^{2} + 16k^{2}\\
  E_{3,k} &: y^{2} = x^{3} - 27x^{2} - 1728k.
\end{align*}
The corresponding morphisms are given by $\phi_{1} : C_{k} \to E_{1,k}$, $\phi_{2} : C_{k} \to E_{2,k}$ and $\phi_{3} : C_{k} \to E_{3,k}$ with
\tiny
\begin{align*}
  \phi_{1}(x : y : z) &= (-xz-yz : xy : z^2)\\
  \phi_{2}(x : y : z) &= (-4x^{3}z^{2} - 4xyz^{3}, -8x^{4}y + 16x^{3}z^{2} + 8y^{3}z^{2} - 12kz^{5} : z^{5})\\
  \phi_{3}(x : y : z) &= \left((x-y)z(16x^{2}y^{2} + 12(x^{3} + x^{2} y + xy^{2} + y^{3})z + 36(x^{2}-xy+y^{2})z^{2}) : \right.\\
&72x^{4}yz + 108x^{4}z^{2} + 64x^{3}y^{3} + 72x^{3}y^{2}z + 108x^{3}z^{3} + 72x^{2}y^{3}z + 216x^{2}y^{2}z^{2} + 72xy^{4}z + 108y^{4}z^{2} + 108y^{3}z^{3} :\\ 
&\left.(x-y)^{3} z^{3}\right).
\end{align*}
\normalsize
For more detail about how these maps were found, see Section~\ref{mapsec}. In the event that one of these three curves (say $E_{i,k}$) has rank zero,
we can determine $C_{k}(\Q)$ by computing $E_{i,k}(\Q)$ followed by
$C_{k}(\Q) = \phi_{i}^{-1}(E_{i,k}(\Q))$. The original motivation for considering this family of curves came from a \href{https://mathoverflow.net/questions/362431}{question on MathOverflow} from Bogdan Grechuk. The question was whether are integral solutions to $x^{3} + x^{2} y^{2} + y^{3} = 7$. The answer is no since the map $\phi_{1} : C_{k} \to E_{1,k}$ sends integral points to integral points. The integral points on $E_{1,7}$ can be found using Magma and one can then verify that none of the preimages on $C_{7}$ of these integral points are themselves integral.

Most of the time if $C_{k}(\Q)$ contains a ``non-trivial'' point, its image on each of $E_{1,k}$, $E_{2,k}$ and $E_{3,k}$ will have infinite order. The main result of this paper is that under the assumption that $k \in \Z$ and one of the $E_{i,k}$ has rank zero,
the only rational points are the obvious ones.

\begin{thm}
\label{main}
Let $k \ne 0$ be an integer and suppose that at least one of $E_{1,k}$, $E_{2,k}$ and $E_{3,k}$ has rank zero over $\mathbb{Q}$. Then
\begin{itemize}
\item If $k \ne 2a^{3} + a^{4}$ for some integer $a$, then $C_{k}(\Q) = \{ (1 : 0 : 0), (0 : 1 : 0) \}$.
\item If $k = 2a^{3} + a^{4}$ for some integer $a \not \in \{ -1, 3 \}$, then $C_{k}(\Q) = \{ (1 : 0 : 0), (0 : 1 : 0), (a : a : 1) \}$.
\item If $k = -1 = 2 \cdot (-1)^{3} + (-1)^{4}$, then $C_{k}(\Q) = \{ (1 : 0 : 0), (0 : 1 : 0), (-1 : 0 : 1), (0 : -1 : 1), (-1 : -1 : 1) \}$.
\item If $k = 135 = 2 \cdot 3^{3} + 3^{4}$, then $C_{k}(\Q) = \{ (1 : 0 : 0), (0 : 1 : 0), (3 : 3 : 1), (-6 : 3 : 1), (3 : -6 : 1) \}$.
\end{itemize}
\end{thm}

\begin{rem}
We are not able to extend our main result to handle the case that $k \in \mathbb{Q}$. One initial roadblock was provably finding all rational points on the curve $y^{2} = x^{6} + 6x^{5} + 39x^{4} + 52x^{3} + 39x^{2} + 6x + 1$, which was handled by Francesca Bianchi and Oana Padurariu using quadratic Chabauty in \cite{BP}.
We were able to rule out the possibility of $E_{1,k}$ or $E_{3,k}$ having a rational point of order $8$,
but we were not able to rule out the possibility of a rational $k \not\in \{ 864, 297/256 \}$ for which $E_{1,k}$ has a rational point of order $5$. For more detail about cases that arise for rational $k$, see Section~\ref{rationalk}.
\end{rem}

The table below lists the integer values of $k$ between $-10$ and $10$, the ranks of $E_{1,k}$, $E_{2,k}$ and $E_{3,k}$ and points $(x : y : 1) \in C_{k}(\Q)$.

\tiny
\begin{tabular}{c|cccc}
$k$ & $\rank E_{1,k}(\Q)$ & $\rank E_{2,k}(\Q)$ & $\rank E_{3,k}(\Q)$ & Some points in $C_{k}(\Q)$\\
\hline
$-10$ & $1$ & $1$ & $1$ & $(1/2 : -9/4 : 1)$, $(-9/4 : 1/2 : 1)$\\
$-9$ & $0$ & $1$ & $0$ & \\
$-8$ & $1$ & $1$ & $1$ & $(-4 : -2 : 1)$, $(0 : -2 : 1)$, $(-2 : 4 : 1)$, $(-2 : 0 : 1)$\\
$-7$ & $0$ & $1$ & $0$ & \\
$-6$ & $1$ & $1$ & $1$ & \\
$-5$ & $1$ & $2$ & $1$ & $(-2 : -1 : 1)$, $(-1 : -2 : 1)$\\
$-4$ & $0$ & $1$ & $1$ & \\
$-3$ & $1$ & $1$ & $1$ & $(-2 : 1 : 1)$, $(-1 : -2 : 1)$\\
$-2$ & $0$ & $0$ & $0$ & \\
$-1$ & $1$ & $0$ & $0$ & $(0 : -1 : 1)$, $(-1 : -1 : 1)$, $(-1 : 0 : 1)$\\
$1$ & $1$ & $1$ & $1$ & $(-3 : -2 : 1)$, $(1 : 0 : 1)$, $(0 : 1 : 1)$, $(1 : -1 : 1)$, $(-1 : 1 : 1)$, $(-2 : -3 : 1)$\\
$2$ & $0$ & $1$ & $0$ & \\
$3$ & $1$ & $1$ & $0$ & $(1 : 1 : 1)$\\
$4$ & $1$ & $1$ & $0$ & \\
$5$ & $1$ & $1$ & $1$ & \\
$6$ & $1$ & $2$ & $1$ & $(-1/2 : 7/4 : 1)$, $(7/4 : -1/2 : 1)$\\
$7$ & $1$ & $1$ & $1$ & \\
$8$ & $1$ & $1$ & $1$ & $(2 : 0 : 1)$, $(0 : 2 : 1)$, $(-4 : 2 : 1)$, $(2 : -4 : 1)$\\
$9$ & $2$ & $1$ & $0$ & \\
$10$ & $0$ & $2$ & $0$ & \\
\end{tabular}

\normalsize
Of the $200$ nonzero integer values of $k$ between $-100$ and $100$, in $99$ cases one of the elliptic curves $E_{i,k}$ has rank zero, and hence Theorem~\ref{main} applies. For the values of $k$ with $1 \leq |k| \leq 1000$, the rank of some $E_{i,k}$ is zero in $961$ out of $2000$ cases. 

The classical method of Chabauty and Coleman applies when the rank of
the Jacobian is less than the genus. If this condition applies to $C_{k}$,
then one of the $E_{i,k}$ must have rank zero and the main theorem applies.
The method of quadratic Chabauty, a special case of the Chabauty-Kim method, proves in an effective way that $C(\Q)$
is finite in other cases. For more detail, see \cite{QC1} and \cite{QC2}. In particular, by Lemma 3.2 of \cite{QC1}, this method will work provided ${\rm rank}~{\rm Jac}(C)(\Q) <
g + \rho - 1$, where $\rho$ is the Picard number of ${\rm
  Jac}(C)/\Q$. Since ${\rm Jac}(C_{k})$ is the product of three
elliptic curves, $\rho \geq 3$, and quadratic Chabauty should be able
to provably determine $C_{k}(\Q)$ provided that the sum of the ranks
of $E_{1,k}$, $E_{2,k}$ and $E_{3,k}$ is $\leq 4$. Nonetheless,
we do not pursue this method since its application to genus $3$
curves is still quite difficult (see \cite{QC3}).

\begin{ack}
The authors used Magma \cite{Magma} V2.26-11 for computations. In addition, we thank Michael Stoll for helpful comments on an earlier draft of this paper.
\end{ack}

\section{Background}
\label{back}
For an elliptic curve $E/\Q$ in the form
$y^{2} + a_{1} xy + a_{3} y = x^{3} + a_{2} x^{2} + a_{4} x + a_{6}$, the set
$E(\Q)$ has the structure of an abelian group. More precisely,
if $P$ and $Q$ are points in $E(\Q)$, let $R = (x,y)$ be the third intersection point of the line through $P$ and $Q$ with the curve, and define
$P+Q = (x,-y-a_{1} x - a_{3})$. The Mordell-Weil theorem states that $E(\Q)$
is a finitely generated abelian group, and Mazur \cite{Mazur} showed that
the torsion subgroup of $E(\Q)$, denoted $E(\Q)_{{\rm tors}}$, must be cyclic
of order $1$, $2$, $3$, $4$, $5$, $6$, $7$, $8$, $9$, $10$ or $12$,
or isomorphic to $\Z/2\Z \times \Z/2k \Z$ for $k \in \{ 1, 2, 3, 4 \}$.

For an elliptic curve $E : y^{2} = x^{3} + Ax + B$ in short
Weierstarss form, we define as usual $j(E) = \frac{6912A^{3}}{4A^{3} +
  27B^{2}}$. We have that two elliptic curves $E_{1}/\Q$ and
$E_{2}/\Q$ are isomorphic over $\overline{\Q}$ if and only if
$j(E_{1}) = j(E_{2})$. If $E : y^{2} = x^{3} + Ax + B$ is an elliptic
curve, the quadratic twist $E_{D}$ of $E$ is given by $E_{D} : y^{2} = x^{3} + AD^{2} x + BD^{3}$.

To determine the torsion subgroup of $E_{1,k}$, $E_{2,k}$ and $E_{3,k}$, we use well-known parametrizations of families of elliptic curves with points of various orders. (See for example page 214 of \cite{Kubert}.)

\begin{lem}
\label{jinvar}
Let $E/\Q$ be an elliptic curve with $j(E) \ne 0, 1728$. Let $f_{i}(t)$ be given in the table below. For
$2 \leq i \leq 7$, if $E/\Q$ has a rational point of order $i$, then there is a rational number $t$ so that
$j(E) = f_{i}(t)$. Conversely, if $j(E) = f_{i}(t)$, then some quadratic twist of $E$ has a rational point of order $i$.

\begin{center}
\begin{tabular}{c|c}
$i$ & $f_{i}(t)$\\
\hline
$2$ & $\frac{t^{3}}{t+16}$\\
\rule{0pt}{3ex}$3$ & $\frac{t^{3}(t+24)}{t-3}$\\
\rule{0pt}{3ex}$4$ & $\frac{-256(t^{2}-t+1/16)^{3}}{t^{4}(t-1/16)}$\\
\rule{0pt}{3ex}$5$ & $\frac{-(t^{4}+12t^{3}+14t^{2}-12t+1)^{3}}{t^{5}(t^{2}+11t-1)}$\\
\rule{0pt}{3ex}$6$ & $\frac{t^{3}(t^{3}-24t-48)^{3}}{(t-6)(t+3)(t+2)^{3}}$\\
\rule{0pt}{3ex}$7$ & $\frac{(t^{2}+t+1)^{3} (t^{6}-5t^{5}-10t^{4}+15t^{3}+30t^{2}+11t+1)^{3}}{t^{7}(t+1)^{7}(t^{3}-5t^{2}-8t-1)}$
\end{tabular}
\end{center}
\end{lem}

We will use the above lemma to determine the integer values of $k$ for which
$j(E_{1,k}) = j(E_{3,k})$ have the above forms. This leads to the equation of a curve (involving the parameters $t$ and $k$). Frequently, this curve has genus $0$ and is isomorphic to $\P^{1}$. This leads to a parametrization $k = f(x)$, $t = g(x)$, and we need to know the rational values of $x$ for which $f(x) \in \Z$. We handle this question using \cite[Lemma~3', p. 72]{SilvermanTate}.

\begin{lem}
\label{stlemma}
Let $\phi(X)$ and $\psi(X)$ be polynomials with integer coefficients and no
common (complex) roots. Let $d$ be the maximum of the degrees of $\phi$ and $\psi$. There is a number $R$, depending on $\phi$ and $\psi$, so that for all rational numbers $\frac{m}{n}$,
\[
  \gcd\left(n^{d} \phi\left(\frac{m}{n}\right), n^{d} \psi\left(\frac{m}{n}\right)\right)
\]
divides $R$. Moreover, if $a_{0}$ is the leading coefficient of $\phi$,
$F(X)$ and $G(X)$ are polynomials with rational coefficients
so that $F(X) \phi(X) + G(X) \psi(X) = 1$ and $A$ is the least common multiples of the denominators of the coefficients of $F$ and $G$, then we can take $R = Aa_{0}^{\deg \phi + \deg \psi}$. 
\end{lem}

In some cases, the above lemma reduces the problem of determining integer
values of a rational function to solving a Thue equation: an equation of
the form $F(x,y) = m$, where $F(x,y)$ is a homogeneous polynomial of degree $\geq 3$. Based on the work of Bilu and Hanrot \cite{BH}, Thue equations can be solved effectively, and we use the Magma \cite{Magma} computer algebra system to solve the ones that arise. Elementary methods can handle cases where $F(x,y)$ is reducible.

\section{Maps from $C_{k}$ to elliptic curves}
\label{mapsec}

For $k \not\in \{0, -\frac{27}{16} \}$, the curve $C_{k} : x^{3} z +
x^{2} y^{2} + y^{3} z = kz^{4}$ is a smooth plane quartic and hence
has genus $3$. The Jacobian $J(C_{k})$ is three-dimensional abelian
variety. As a consequence of the maps given in the introduction, this
Jacobian decomposes (up to isogeny) as the product of three elliptic
curves. We wish to explain briefly how these maps were found.

First, the affine equation for $C_{k}$ is symmetric in $x$ and $y$. So letting
$s = x+y$ and $t = xy$, the equation of $C_{k}$ becomes $s^{3} - 3st + t^{2} = k$.
Setting $X = -s$ and $Y = t$ the equation becomes $Y^{2} + 3XY = X^{3} + k$
and the map $(x,y) \mapsto (s,t)$ is a degree $2$ map.

In addition, the curve $C_{k}$ has an automorphism of order $3$ defined
by $\phi(x : y : z) = (\zeta_{3} x : \zeta_{3}^{-1} y : z)$ where $\zeta_{3} = e^{2 \pi i / 3}$. The quotient of $C_{k}$ by this subgroup of ${\rm Aut}(C_{k})$ generated by $\phi$ is a genus $1$ curve which is isomorphic to $E_{2,k}$. 

Since $J(C_{k})$ has maps to at least two non-isomorphic elliptic curves,
Poincar\'e's complete reducibility theorem (see \cite[Proposition 12.1]{Milne}) implies that $J(C_{k})$ is
isogenous to the product of three elliptic curves. To find the third
elliptic curve, we choose several small values of $k$ and compute the
numerator of the zeta function of $C_{k}/\F_{p}$ for several lots of
small primes $p$ using the algorithm of Tuitman (see \cite{Tuitman1}
and \cite{Tuitman2}) and implemented in Magma. We observe that for $p
\equiv 1 \pmod{3}$, the numerator of the zeta function has a repeated quadratic
factor, while for $p \equiv 2 \pmod{3}$ it usually does not have a repeated quadratic factor. This suggests that $J(C_{k})/\F_{p}$ has a repeated factor when
$p \equiv 1 \pmod{3}$ and suggests the remaining elliptic curve is the quadratic
twist by $-3$ of either $E_{1,k}$ or $E_{2,k}$. Checks for a few values of $k$
indicate that the new elliptic curve is the quadratic twist of $E_{1,k}$ by $-3$.

We proceed to search for a map from $C_{k}$ to $E_{3,k}$. We choose
two small positive integers $a$ and $b$ and let $k = a^{3} + a^{2}
b^{2} + b^{3}$. This forces the point $(a : b : 1)$ to be on
$C_{k}$. For values of $k$ that have this form, we test if $E_{3,k}$
has rank $1$ and trivial torsion subgroup. If it does, we search for
points of small height on $E_{3,k}$ in the hope that the simplest
point in $E_{3,k}(\Q)$ is the image of $(a : b : 1)$ under the map $C_{k} \to E_{3,k}$. This leads easily to the map given in the introduction. Once the map is found, it is straightforward to verify that it is correct, and
that $\phi_{3} : C_{k} \to E_{3,k}$ has degree $6$.

\section{Proof of Theorem~\ref{main}}
\label{proofsec}

Suppose that $(x : y : 1)$ is a rational point on $C_{k}$. It follows that $\phi_{1}(x : y : 1) \in E_{1,k}(\Q)$, $\phi_{2}(x : y : 1) \in E_{2,k}(\Q)$, and $\phi_{3}(x : y : 1) \in E_{3,k}(\Q)$. Either all three points
have infinite order, or one or more of points are torsion points. The
proof of Theorem~\ref{main} follows from classifying the values of $k$
for which these curves have torsion points, together with a finer analysis
of rational points of order $2$ on $E_{1,k}$ and $E_{3,k}$.

First, we note that $E_{2,k}(\Q)$ has positive rank unless $k = -1$ or
$k = -2$.  In particular, $E_{2,k} : y^{2} = x^{3} + 4x^{2} + 16k^{2}$
and so $P = (0,4k)$ is a rational point on $E_{2,k}$. Computing the
coordinates of $2P$, $3P$, $\ldots$, $12P$ as rational functions in
$k$ we find that $P$ can have finite order if and only if $k = -1$ or
$k = -2$. When $k = -1$, $P$ has order $5$ and when $k = -2$, $P$ has order $6$.
It suffices then to determine when $\phi_{1}(x : y : 1)$ and $\phi_{3}(x : y : 1)$ have finite order.

If either of these points has finite order, its order must be in $\{ 1, 2, 3, 4, 5, 6, 7, 8, 9, 10, 12 \}$. However, we will show that there are only finitely many integers $k$ for which $E_{1,k}$ or $E_{3,k}$ has a rational point of order
$3$, $4$, $5$ or $7$, and this will suffice.

{\bf Case 1}: Either $\phi_{1}(x : y : 1)$ or $\phi_{3}(x : y : 1)$ has order $1$.

The point of order 1 should be the point $(0:1:0)$ in projective space. For $E_{1,k}$, the corresponding morphism is $\phi_{1}(x : y : 1) = (-x-y : xy : 1)$. This cannot equal $(0 : 1 : 0)$. So only $E_{3,k}$ will be in consideration.

For $E_{3,k}$, the corresponding morphism is
\small
\begin{align*}
\phi_{3}(x:y:1) &=\left.(16x^{2}y^{2} + 12(x^{3} + x^{2} y + xy^{2} + y^{3}) + 36(x^{2}-xy+y^{2}))(x-y): \right.\\
&\left.72x^{4}y + 108x^{4} + 64x^{3}y^{3} + 72x^{3}y^{2} + 108x^{3} + 72x^{2}y^{3} + 216x^{2}y^{2} +\right.\\
&\left.72xy^{4} + 108y^{4} + 108y^{3} : (x-y)^{3})\right..\\
\end{align*}
\normalsize
This equals $(0 : 1 : 0)$ if and only if $x = y$. This implies that $(x : y : 1)$ on $C_{k}$ satisfies $2x^{3} + x^{4} = k$.

{\bf Case 2}: Either $\phi_{1}(x : y : 1)$ or $\phi_{3}(x : y : 1)$ has order $2$.

We will show that this only occurs for $k \in \{ -17, -1, 135, 368 \}$.
If $k \in \{ -17, 368 \}$, then $E_{1,k}$, $E_{2,k}$ and $E_{3,k}$ all have positive rank. 

If $k = -1$, we find that $E_{3,k}(\Q) \cong \Z/2\Z$ and the two points
are $(0 : 1 : 0)$ and $(24 : 0 : 1)$. We compute all rational points on $C_{-1}$ which are preimages of these two points under $\phi_{3}$. 
For $(0 : 1 : 0)$ we find $(-1 : -1 : 1)$, $(0 : 1 : 0)$ and $(1 : 0 : 0)$. For $(24 : 0 : 1)$ we find $(-1 : 0 : 1)$ and $(0 : -1 : 1)$. These are all the rational points on $C_{-1}$.

If $k = 135$, we also have $E_{3,k}(\Q) \cong \Z/2\Z$ and the rational points are $(0 : 1 : 0)$ and $(72 : 0 : 1)$. The rational points on $C_{135}$ that are preimages of $(0 : 1 : 0)$ are $(3 : 3 : 1), (0 : 1 : 0), (1 : 0 : 0)$ and the rational points on $C_{135}$ that are preimages of $(72 : 0 : 1)$ are $(-6 : 3 : 1)$ and $(3 : -6 : 1)$.

First consider the case of $E_{1,k}$. To easily identify points of order $2$,
we rewrite the model as $E_{1,k}' : Y^{2}=X^{3}-\frac{27}{16}X+(k+\frac{23}{32})$.
The map $g : E_{1,k} \to E_{1,k}'$ is given by
$g(x,y) = (x + \frac{3}{4}, y + \frac{3}{2}x)$. Now we want to determine the values of $k \in \mathbb{Z}$ for which there is some $P \in C_{k}(\mathbb{Q})$ so that the $y$-coordinate of $g \circ \phi_{1}(P) = 0$.

If $P = (x : y : 1)$ we obtain the equation $xy = \frac{3}{2}(x+y)$ and we
also have $k = x^{3} + x^{2}y^{2} + y^{3} \in \Z$. We can express $x$ in terms of $y$, and find $x=\frac{3y}{2y-3}$. The assumption that $k \in \Z$ implies $\left(\frac{3y}{2y-3}\right)^{3}+\left(\frac{3y}{2y-3}\right)^{2}y^{2}+y^{3}$ is an integer. Write $y=\frac{a}{b}$ where $a,b \in \mathbb{Z}$ with $\gcd(a,b)=1$ and $b > 0$. Thus $\frac{27a^{3}b^{3}+9a^{4}b(2a-3b)+a^{3}(2a-3b)^{3}}{b^{3}(2a-3b)^{3}}$ is an integer. So $b \mid 27a^{3}b^{3}+9a^{4}b(2a-3b)+a^{3}(2a-3b)^{3}$ which implies $b \mid a^{3}(2a-3b)^{3}$. Since $\gcd(a,b) = 1$, this implies that $b \mid 8$ and so
$b \in \{ 1, 2, 4, 8 \}$.

If $b = 1$, then $2a-3 \mid 27a^{3}$ and this implies $2a-3 \mid 8 \cdot 27a^{3}$. Since $2a-3 \mid (2a)^{3}-3^{3}$ and thus $2a - 3 \mid 27 ((2a)^{3} - 3^{3})$. It follows that $2a - 3 \mid 729$. There are only finitely many $a$ with this property and for these values of $a$ we find $k \in \{ -17, 0, 135, 368 \}$.

If $b > 1$, then $b$ is even and $2a - 3b \mid (2a)^{3} - (3b)^{3}$ implies that
$2a - 3b \mid (\frac{b}{2})^{3}(2a)^{3}-(\frac{b}{2})^{3}(3b)^{3}$. Also, the assumption that $k \in \Z$ implies that
$2a-3b \mid 27a^{3}b^{3}$. Since $2a - 3b \mid 27(\frac{b}{2})^{3}(2a)^{3}-27(\frac{b}{2})^{3}(3b)^{3}$, we get that $2a-3b \mid \frac{729}{8}b^{3}$. Since $b \in \{ 2, 4, 8 \}$ there are only finitely many possibilities for $a$. We find these using Magma and again determine that $k \in \{ -17, 0, 135, 368 \}$.

Now consider the case of $E_{3,k}$. If $\phi_{3} : C_{k} \to E_{3,k}$
is the morphism given above and the image of $P = (x : y : 1)$ on $E_{3,k}$
has order $2$, then the numerator of the $y$-coordinate of $\phi_{3}(x : y : 1)$ must be zero. This equation defines a curve 
\begin{align*}
  C &: 72x^{4}yz + 108x^{4}z^{2} + 64x^{3}y^{3} + 72x^{3}y^{2}z + 108x^{3}z^{3}\\ 
  &+ 72x^{2}y^{3}z + 216x^{2}y^{2}z^{2} + 72xy^{4}z + 108y^{4}z^{2} + 108y^{3}z^{3} = 0
\end{align*}
of genus $2$. This curve $C$ is isomorphic to $y^{2} = x^{6} + 6x^{5} + 39x^{4} + 52x^{3} + 39x^{2} + 6x + 1$.
In \cite{BP}, Bianchi and Padurariu use quadratic Chabauty to show that this curve has precisely 12 rational points. These points correspond to $k = 0$, $\infty$, $1$ and $135$.
Bianchi and Padurariu use the isomorphic model $y^{2} = x^{6} + 3x^{4} - 9x^{2} + 9$ and the complete list of rational points on this curve is given in the list {\tt pts[354]} available in the file \href{https://github.com/oana-adascalitei/MWSieveForDatabase/blob/main/allcurvesOutput.m}{{\tt allcurvesOutput.m}} at their \href{https://github.com/oana-adascalitei/MWSieveForDatabase}{Github repository}. We are able to give an independent classification of the integer values of $k$ for which $E_{3,k}$ has a rational point of order $2$ by relating the problem to the rational points on the elliptic curve $y^{2} = x^{3} - 12x + 20$ for which a particular function $f : E \to \P^{1}$ takes integer values. The finite set of such points can be determined using Corollary IX.3.2.2 of \cite{Silverman} and the determination of the $S$-integral points on $y^{2} = x^{3} - 12x + 20$ for $S = \{ 2 \}$ using Magma.

\begin{rem}
  There are infinitely many integers $k$ for which $E_{1,k}$ and $E_{3,k}$
  have rational points of order $2$. This is the reason for considering
  the more delicate question of when there is a rational point $P \in C_{k}(\Q)$ whose image on $E_{1,k}$ or $E_{3,k}$ has order $2$.
\end{rem}

{\bf Case 3}: $E_{1,k}(\Q)$ or $E_{3,k}(\Q)$ contains a rational point of order $3$.

We will show that this occurs only for $k \in \{ -72, -2, 2058, -56000 \}$. If $k = -72$,
then $E_{3,k}(\Q) = \{ (0 : 1 : 0) \}$ and pulling back this point
we find that $C_{k}(\Q) = \{ (0 : 1 : 0), (1 : 0 : 0) \}$. If $k = -2$,
we find that $E_{1,k}(\Q) = \{ (0 : 1 : 0) \}$ and pulling back this point
we find that $C_{k}(\Q) = \{ (0 : 1 : 0), (1 : 0 : 0) \}$.
If $k = 2058$, $E_{1,k}$, $E_{2,k}$ and $E_{3,k}$ all have rank $1$. If $k = -56000$, $E_{3,k}(\Q) = \{ (0 : 1 : 0) \}$ and pulling back we find that $C_{k}(\Q) = \{ (0 : 1 : 0), (1 : 0 : 0) \}$. 

If $E_{1,k}$ or $E_{3,k}$ has a rational point of order $3$, then there is a rational number $t$ so that
\[
j(E_{1,k}) = j(E_{3,k}) = \frac{-3^{7}}{16 k^{2} + 27k} = \frac{t^{3} (t+24)}{t-3}.
\]
The equation above in $k$ and $t$ parametrizes a curve of genus zero with a non-singular rational point. We parametrize this curve and find that
\[
  k = \frac{u^{3} (-2u-3v)}{(u+v)^{3} (u-3v)}.
\]  
for some integers $u$ and $v$ with $\gcd(u,v) = 1$. 

Now, we apply Lemma~\ref{stlemma}. Let $\psi(x) = -2x^{4} - 3x^{3}$
and $\phi(x) = (x+1)^{3} (x-3)$. We find there exist $F(x)$ and $G(x)$ that $F(x) \phi(x) + G(x) \psi(x) = 1$ and obtain that
\[
  \gcd(v^{4} \phi(u/v), v^{4} \psi(u/v))
\]
divides $243$.

Since we assume that $k$ is an integer, it follows that $(u+v)^{3} (u-3v)$ divides the numerator, and hence $(u+v)^{3} (u-3v)$ divides $243$. 

For each divisor
$d$ of $243$, we solve the equation $(u+v)^{3} (u-3v) = d$ by enumerating divisors $d_{1}$ of $d$ that are cubes
and solving the system $(u+v)^{3} = d_{1}$, $u-3v = d/d_{1}$. Using this process, we find that the only integral values of $k$ of the form
\[
  \frac{u^{3} (-2u-3v)}{(u+v)^{3} (u-3v)}
\]
are $k = -72$, $k = -2$, $k = 0$, $k = 2058$, and $k = -56000$. For these values of $k$
we compute the torsion subgroup of $E_{1,k}$ and $E_{3,k}$ and obtain the desired
result.

{\bf Case 4}: $E_{1,k}(\Q)$ or $E_{3,k}(\Q)$ contains a rational point of order $4$.

We will show that this only occurs if $k = 135$, and we previously found $C_{135}(\Q)$ in Case 2.

An elliptic curve has a point of order $4$ if and only if it can be put in the form $y^{2}+xy+ty = x^{3}+ty^{2}$ for some $t \in \mathbb{Q}$. So if $k$ is an integer and $E_{1,k}$ or $E_{3,k}$ has a point of order $4$, then
\[
j(E_{1,k}) = j(E_{3,k}) = \frac{-3^{9}}{16k^{2}+27k} = \frac{-256(t^{2}-t+1/16)^{3}}{t^{4}(t-1/16)}
\]
for some rational number $t$.

This equation in $k$ and $t$ defines a genus zero curve which is birational to $\P^{1}$ and it follows that $k$ must have the form
\[
k=\frac{-\frac{27}{16}t^{6}+\frac{243}{16}t^{4}}{(t^{2}-3)^{3}}.
\]
Let $\mathbb{Q}(t)$ be the collection of rational functions. We wish to
determine the $x$-coordinates of points of $E_{1,k}$ and $E_{3,k}$ with order $4$
in $\Q(t)$. We will use the Magma command {\tt DivisionPolynomial(E,4)} to compute these and find that the roots for $E_{3,k}$ are $(18t^{2} - 54t)/(t^{2} - 3)$ and $(18t^{2} + 54t)/(t^{2} - 3)$. The roots for $E_{1,k}$ are $(-3/2t^{2} - 9/2t)/(t^{2} - 3)$ and $(-3/2t^{2} +9/2t)/(t^{2} - 3)$.

Write $t=a/b$ for integers $a$ and $b$ with $\gcd(a,b)=1$. We get
\[
k = \frac{-\frac{27a^{6}}{16b^{6}}+\frac{243a^{4}}{16b^{4}}}{(\frac{a^{2}}{b^{2}}-3)^{3}}=\frac{-\frac{27}{16}a^{6}+\frac{243}{16}a^{4}b^{2}}{(a^{2}-3b^{2})^{3}}=\frac{-27a^{6}+243a^{4}b^{2}}{16(a^{2}-3b^{2})^{3}}.
\]
Since $k$ should be an integer, $16(a^{2}-3b^{2})^{3} \mid (-27a^{6}+243a^{4}b^{2})$. Now we consider two situations. 

If $3 \nmid a$, then because $\gcd(a,b)=1$, one can show that $a^{2}-3b^{2} \mid a^{2}-9b^{2}$ implies $a^{2}-3b^{2} \mid 2$. Considering the congruence situation modulo $3$, the only possibilities are $a^{2}-3b^{2}=1$ and $a^{2}-3b^{2}=-2$. For $a^{2}-3b^{2}=-2$, both $a$ and $b$ should be odd numbers. Otherwise, either $a^{2}-3b^{2}$ is odd or $\gcd(a,b)$ is greater than 1. If $a^{2}-3b^{2}$ is odd, then $-128 \mid \left(-27a^{4}(a^{2}-9b^{2})\right)$ and $-128 \mid (-2-6b^{2})$. Then $64|1-3b^{2}$. Because $b$ is an odd number, $1-3b^{2} \equiv 6 \pmod{8}$ and this shows that $64 \nmid 1-3b^{2}$. The conclusion is if $3 \nmid a$, then $a^{2}-3b^{2}=1$.

If $3 \mid a$, then let $a =3c$. Because $\gcd(a,b)=1$, $3 \nmid b$ and $\gcd(b,c)=1$. In a similar way, one can derive $3c^{2}-b^{2}=-1$ or $-2$ and so $a^{2}-3b^{2} =-3$ or $a^{2}-3b^{2}=-6$. This shows that if $3 \mid a$, then $a^{2}-3b^{2} =-3$ or $a^{2}-3b^{2} =6$.

In particular, in any case we have $a^{2} - 3b^{2} \in \{ 1, -3, 6 \}$.

The conditions on $k$ and the rational number $t$ ensures that $E_{1,k}$ and $E_{3,k}$ have points of order $4$ whose $x$-coordinate is rational. The condition that the $y$-coordinate is rational gives rise to an elliptic curve. In the case of $E_{3,k}$, there is a rational number $u$ so that $u^{2} = 3(t-1)(t-3)(t^{2}-3)$ or $u^{2} = 3(t+1)(t+3)(t^{2}-3)$. In the first case the elliptic curve has rank zero and torsion subgroup $\Z/2\Z$. This means that the only solutions are $t = 1$ and $t = 3$. The former gives $k = -27/16$, which is not integral, and the latter gives $k = 0$, and we assume $k \ne 0$. The second case is similar; the elliptic curve has rank zero and torsion subgroup $\Z/2\Z$ and the two points correspond to $k = -27/16$ and $k = 0$ again. So we have shown that there are no integral values of $k$ for which $E_{3,k}$ has a point of order $4$. 

In the case of $E_{1,k}$, the $x$-coordinate of a point of order $4$ is $\frac{-(3/2)t^{2}+(9/2)}{t^{2}-3}$ or $\frac{-(3/2)t^{2}-(9/2)}{t^{2}-3}$. In order for the $y$-coordinate to be rational, we need a rational solution to $u^{2}=-(t+1)(t+3)(t^{2}-3)$, or $u^{2}=-(t-1)(t-3)(t^{2}-3)$ in the other case.
Both cases give rise to the elliptic curve $y^{2}=x^{3}-12x$, which has rank $1$. Also, we can change $u^{2}=-(t-1)(t-3)(t^{2}-3)$ to $u^{2}=-(-t+1)(-t+3)(t^{2}-3)$. Negating $t$ will not affect the value $k$ because the powers of $t$ in the polynomial are even, we only focus on $u^{2}=-(t+1)(t+3)(t^{2}-3)$.

This can be transformed into $y^{2}=x^{3}-12x$ where $x=\frac{2a+6b}{a+b}$. Rather than use the technique of Corollary IX.3.2.2 of \cite{Silverman}, we will use an elementary argument based on Pythagorean triples to find the rational points on $y^{2} = x^{3} - 12x$ for which $k$ is an integer. 

We have that
\[
y^{2}=x^{3}-12x=\frac{-16 (a+b)(a+3b)(a^{2}-3b^{2})}{(a+b)^{4}} = \frac{(-16 \text{ or }  48 \text{ or }-96)(a+b)(a+3b)}{(a+b)^{4}}
\]
is a perfect square using that $a^{2}-3b^{2} = 1$,$-3$, or $6$.

Since it is a perfect square, we can use Pythagorean triples to show that there are only finitely many options for the coprime integers $a$ and $b$.

First, if $a^{2}-3b^{2}=1$, then $\frac{-16(a+b)(a+3b)}{(a+b)^{4}}$ is a perfect square. Write $-(a+b)(a+3b)=u^{2}$ for a positive integer $u$. So $-(a^{2}+4ab+3b^{2})=b^{2}-(a+2b)^{2}=u^{2}$. Then $a+2b,u,b$ is a Pythagorean triple. By Theorem 8.7 on page 101 of \cite{Marshall}, there exist relatively prime positive integers $s$ and $t$, one even and one odd, such that $b = s^{2}+t^{2}$, and
$\{ u, a+2b \} = \{ 2st, s^{2}-t^{2} \}$.

If $a$ is even then $a+2b = 2st$, and so $a = 2st-2s^{2}-2t^{2}$. So $a^{2}-3b^{2} = s^{4}-8s^{3}t+6s^{2}t^{2}-8st^{3} + t^{4} = 1$. By using the Thue function in Magma, we find $(s,t)$ can be $(1,0)$, $(-1,0)$, $(0,1)$, or $(0,-1)$. In this situation, $a=-2$ and $b=1$ and so $t = -2$ and this yields $k = 135$. 

If $a$ is odd then $a+2b = s^{2}-t^{2}$, $a = s^{2}-t^{2} - 2s^{2}-2t^{2}= -s^{2}-3t^{2}$. In this situation, $1 = a^{2}-3b^{2} = -2s^{4} + 6t^{4}$ is an even number, which is a contradiction.

Second, if $a^{2}-3b^{2}=-3$, then $\frac{48(a+b)(a+3b)}{(a+b)^{4}}$ is a perfect square. Then $3(a+b)(a+3b) = 3a^{2}+12b^{2}+9b^{2} = (2a+3b)^{2} -a^{2}$. Let $(2a+3b)^{2} -a^{2} = u^{2}$. Since $\gcd(a,b)=1$, one can also deduce that $\gcd(2a+3b,a) =1$ and $(2a+3b,u,a)$ is a Pythagorean triple. There are relatively prime integers $s$ and $t$ such that $2a+3b=s^{2}+t^{2}$. In this situation $b$ is odd. 

If $a$ is odd, then $a = s^{2}-t^{2}$. Since $2a+3b=s^{2}+t^{2}$, $b=t^{2}-\frac{s^{2}}{3}$. Then $a^{2}-3b^{2}= \frac{2}{3}s^{4}-2t^{4}$ is equal to $-3$. This implies that
$2s^{4} - 6t^{4} = -9$, which is a contradiction since $2s^{4} - 6t^{4}$ is even.

If $a$ is an even number, then $a=2st$, so $b=\frac{s^{2}+t^{2}-4st}{3}$. Then $a^{2}-3b^{2}=4s^{2}t^{2}-\frac{1}{3}(s^{4}-8s^{3}t+18s^{2}t^{2}-8st^{3}+t^{4})=-\frac{1}{3}s^{4}+\frac{8}{3}s^{3}t-2s^{2}t^{2}+\frac{8}{3}st^{3}-\frac{1}{3}t^{4}$ is equal to $-3$. This shows that $s^{4}-8s^{3}t+6s^{2}t^{2}-8st^{3}+t^{4}$ is equal to $9$. By using the Thue function in Magma, we find that there are no solutions in this case.

Third, if $a^{2}-3b^{2}=6$, then $\frac{-96(a+b)(a+3b)}{(a+b)^{4}}$ is a perfect square. One can let $-6(a+b)(a+3b)(a+b)^{4} = d^{2}$. In the previous proof, we showed that both $a$ and $b$ are odd and $\gcd(a,b)=1$. Let $k = \gcd(a+b, a+3b)$. Then $k \mid a+3b - (a+b)$ which means $k \mid 2b$. Also, $k \mid 2(a+b)-2b$ means $k \mid 2a$. This shows that $k \mid 2b$ and $k \mid 2a$ and $\gcd(a,b)$=1. We get $k=2$ and $\gcd(\frac{a+b}{2},\frac{a+3b}{2})=1$. We also get the equation $-6((a+b)/2)((a+3b)/2)=(d/2)^{2}$.

Because $\gcd(\frac{a+b}{2},\frac{a+3b}{2})=1$ and $-6((a+b)/2)((a+3b)/2)=(d/2)^{2}$, $3 \mid \frac{a+3b}{2}$ and there are four possible options:
either $\frac{a+b}{2}$ is a square, the negative of a square, twice a square, or
the negative of twice a square.

We can represent the four situations as $\frac{a+b}{2} = cs^{2}$ and $\frac{a+3b}{2} = -(\frac{6}{c})t^{2}$ with $c \in \{ \pm{1}, \pm{2} \}$. Then $b =\frac{a+3b}{2}-\frac{a+b}{2} =(\frac{6}{c})t^{2}-cs^{2} $ and $a = 3(\frac{a+b}{2})-\frac{a+3b}{2} = 3cs^{2}+(\frac{6}{c})t^{2}$. We also know $a^{2}-3b^{2} = 6$ in this situation. This gives $a^{2}-3b^{2} = 6c^{2}s^{4} - \frac{72}{c^{2}}t^{4}=6$. Since $c \in \{ \pm{1}, \pm{2} \}$ we get $s^{4}-12t^{4}=1$ or $4s^{4}-3t^{4}=1$.

Using the Thue function in Magma, we find that the only solutions
are $(c,a,b) = (1,3,-1)$, $(-1,-3,1)$, $(2,9,-5)$, and $(-2,-9,5)$. These
yield $t = -3$ or $t = -9/5$. These lead to $k = 0$ and $k = 59049/8$.

In conclusion, the only solution we found which yields a value of $k$ which is a nonzero integer was $(a,b) = (-2,1)$.
In this case $k = 135$ and $E_{1,k}$ has a rational point of order 4.

{\bf Case 5}: $E_{1,k}(\Q)$ or $E_{3,k}(\Q)$ contains a rational point of order $5$.

We will show that this only occurs if $k = 864$. We have $864 = 2 \cdot (-6)^{3} + (-6)^{4}$. We have that $E_{3,k}(\Q)$ is trivial, and therefore every rational point $P \in C_{864}(\Q)$ maps to $(0 : 1 : 0)$ on $E_{3,k}$. The calculation in Case 1 shows that $C_{864}(\Q) = \{ (1 : 0 : 0), (0 : 1 : 0), (-6 : -6 : 1) \}$.

If $E_{1,k}$ or $E_{3,k}$ has a rational point of order $5$, then there is a rational number $t$ so that
\[
j(E_{1,k}) = j(E_{3,k}) = \frac{-3^{9}}{16 k^{2} + 27k} = \frac{-(t^{4}+12t^{3}+14t^{2}-12t+1)^3}{t^5(t^2+11t-1)} 
\]
The equation above in $k$ and $t$ parametrizes a curve of genus zero with a non-singular rational point. Parametrizing this curve shows that 
\[
  16k^{2}+27k = \frac{3^9(u^{7}v^{5}+11u^{6}v^{6}-u^{5}v^{7})}{(u^{4}+12u^{3}v+14u^{2}v^{2}-12uv^{3}+v^{4})^{3}}
\]  
for some integers $u$ and $v$ with $\gcd(u,v) = 1$. 

Now, we apply Lemma~\ref{stlemma}. Let $\psi(x) = 3^9 (x^{7}+11x^{6}-x^{5})$
and $\phi(x) = (x^{4}+12x^{3}+14x^{2}-12x-1)^3$. We find that there exist $F(x)$ and $G(x)$ such that $F(x) \phi(x) + G(x) \psi(x) = 1$ and obtain that
\[
  \gcd(v^{12} \phi(u/v), v^{12} \psi(u/v))
\]
divides $2460375$.

Since we assume that $k$ is an integer, it follows that $(u^{4}+12u^{3}v+14u^{2}v^{2}-12uv^{3}+v^{4})^{3}$ divides the numerator, and hence $(u^{4}+12u^{3}v+14u^{2}v^{2}-12uv^{3}+v^{4})^{3}$ divides $2460375$. For each divisor
$d$ of $2460375$ which is a cube, we use Magma to solve the Thue equation $(u^{4}+12u^{3}v+14u^{2}v^{2}-12uv^{3}+v^{4})^{3} = d$.
Using this process, we find that the only integral values of the form
\[
  16k^{2} + 27k = \frac{3^9(u^{7}v^{5}+11u^{6}v^{6}-u^{5}v^{7})}{(u^{4}+12u^{3}v+14u^{2}v^{2}-12uv^{3}+v^{4})^{3}}
\]
are $0$ and $11967264$. By solving the quadratic equation $16k^{2}+27k = 11967264$, we find $k=13581/16$ and $k=864$. The only integral value is therefore $k = 864$ and for this value of $k$
we find that $E_{1,k}$ has a rational point of order $5$ and $E_{3,k}$ does not.

{\bf Case 6}: $E_{1,k}(\Q)$ or $E_{3,k}(\Q)$ contains a rational point of order $7$.

We will show that this does not occur, even if $k$ is a rational number.

If $E_{1,k}$ or $E_{3,k}$ has a rational point of order $7$, then there is a rational number $t$ so that
\[
j(E_{1,k}) = j(E_{3,k}) = \frac{-3^{9}}{16 k^{2} + 27k} = \frac{(t^{2}+t+1)^{3} (t^{6} - 5t^{5} - 10t^{4} + 15t^{3} + 30t^{2} + 11t + 1)^{3}}{t^{7} (t+1)^{7} (t^{3} - 5t^{2} - 8t - 1)}. 
\]
This implies that $k^{2} = -\frac{3^{9}}{j} + 729/64$ must be a square, and this yields that
$t^{8} - 4t^{7} - 14t^{6} + 35t^{4} + 56t^{3} + 42t^{2} + 12t + 1$ must be a square. Let
\[
a = \sqrt{t^{8} - 4t^{7} - 14t^{6} + 35t^{4} + 56t^{3} + 42t^{2} + 12t + 1}.
\]
Computing the $7$-division polynomial for $E_{1,k}$ over $\Q(a)$ shows that in order for $E_{1,k}$ to have a rational point of order $7$, we must have that $a$ is a square (and for $E_{3,k}$ to have a rational point of order $7$, we must have that $-3a$ is a square).
Letting $y^{2} = a$, we therefore seek to know the rational points on $y^{4} = t^{8} - 4t^{7} - 14t^{6} + 35t^{4} + 56t^{3} + 42t^{2} + 12t + 1$. The automorphism group of this curve over $\Q$
is isomorphic to $\Z/6\Z$ and the quotient by the automorphism group of order $3$ is the genus $3$ curve $X : y^{2} = x^{8} - 6x^{7} - 7x^{6} - 42x^{5} - 42x^{3} -7x^{2} -6x - 1$. This curve has four automorphisms over $\Q$, and one of the quotients is the genus $2$ curve $y^{2} = x^{6} - 6x^{5} - 15x^{4} + 60x^{2} + 96x - 64$. 
The Jacobian of this curve has rank $1$ and Magma's implementation of Chabauty shows that there are precisely four rational points on this genus $2$ curve. Pulling these points back to $X$ shows that the only rational points have $t = \infty$, $t = 0$ or $t = -1$, and these correspond to $k = 0$ and $k = \infty$.
On the other hand, the curve $9y^{4} = t^{8} - 4t^{7} - 14t^{6} + 35t^{4} + 56t^{3} + 42t^{2} + 12t + 1$ has no $\Q_{3}$-points. As a consequence, there are no values of $k$ for which $E_{1,k}$ or $E_{3,k}$ has a rational point of order $7$.

We conclude with a table below of the special cases of $k$ that arose in the proof above.

\scriptsize
\begin{tabular}{l|lllll}
$k$ & Is $k = 2a^{3} + a^{4}$? & $E_{1,k}(\Q)$ & $E_{2,k}(\Q)$ & $E_{3,k}(\Q)$ & Some points $(x : y : 1) \in C_{k}(\Q)$\\
\hline
$-1$ & Yes, $a = -1$ & $\Z/2\Z \times \Z$ & $\Z/5\Z$ & $\Z/2\Z$ & $\{ (0 : -1 : 1), (-1 : -1 : 1), (-1 : 0 : 1) \}$\\
$-2$ & No & Trivial & $\Z/6\Z$ & $\Z/3\Z$ & None.\\
$-17$ & No & $\Z/2\Z \times \Z$ & $\Z \times \Z$ & $\Z/2\Z \times \Z$ & $\{ (2 : -5 : 1), (1 : -3 : 1), (-5 : 2 : 1), (-3 : 1 : 1) \}$\\
$-72$ & No & $\Z/3\Z$ & $\Z$ & Trivial & None.\\
$135$ & Yes, $a = 3$ & $\Z/4\Z$ & $\Z$ & $\Z/2\Z$ & $\{ (3 : -6 : 1), (3 : 3 : 1), (-6 : 3 : 1) \}$\\
$368$ & No & $\Z/2\Z \times \Z$ & $\Z \times \Z$ & $\Z/2\Z \times \Z \times \Z$ & $\{ (6 : 2 : 1) : (2 : 6 : 1) \}$\\
$864$ & Yes, $a = -6$ & $\Z/5\Z$ & $\Z$ & Trivial & $\{ (-6 : -6 : 1) \}$.\\
$2058$ & No & $\Z$ & $\Z/3\Z \times \Z$ & $\Z/3\Z \times \Z$ & None.\\
$-56000$ & No & $\Z/3\Z \times \Z$ & $\Z \times \Z$ & Trivial & None.\\
\end{tabular}
\normalsize

\section{Extending the main result to $k \in \mathbb{Q}$}
\label{rationalk}

In this section, we describe the challenges to extending Theorem~\ref{main}
to apply to all $k \in \mathbb{Q}$. If $k = 0$ or $k = -27/16$, then $C_{k}$ is singular and we disregard these cases.

For each positive integer $n$, there are curves $C_{1,n}$ and $C_{3,n}$
that parametrize points $(a : b : 1)$ on some $C_{k}(\Q)$ whose image
on $E_{1,k}$ has order $n$ (respectively, for which the image on $E_{3,k}$ has
order $n$). Each such curve has an involution swapping $a$ and $b$.

In addition, for each positive integer $n$, there are curves $D_{1,n}$
and $D_{3,n}$ that parametrizes rational values of $k$ for which $E_{1,k}(\Q)$
(respectively $E_{3,k}(\Q)$) contains a rational point of order $n$.

To extend Theorem~\ref{main} it is necessary to provably find the rational points on these curves for $n \in \{ 2, 3, 4, 5, 6, 7, 8, 9, 10, 12 \}$. We document here our attempts to find such points.

The curve $C_{1,2}$ has affine model $8ab - 12a - 12b = 0$. This curve is
isomorphic to $\mathbb{P}^{1}$, and so if
$k = \frac{(27/8) t^{4} (t^{2} - 3/2 t + 3/2)}{(t-1)^{3}}$,
there is a rational point on $C_{k}$, namely
\[
  \left(\frac{3t}{2t-2} : \frac{3t^{2} - 3t}{2t - 2} : 1\right)
\]
whose image on $E_{1,k}$ has order $2$.

The curve $C_{1,3}$ has genus $1$ and is isomorphic to $y^{2} + xy + y = x^{3} - x$, which has rank zero. The only rational points on $C_{1,3}$ correspond
to $k = 0$ or $k = \infty$.

The curve $C_{1,4}$ has genus $4$ and has a number of rational points corresponding to $k = \infty$, $k = 0$, $k = -27/16$ and $k = 135$ on it. We have not been able to provably determine its rational points. The quotient by the involution on $C_{1,4}$ is the elliptic curve $y^{2} = x^{3} - 12x$ and a model for $C_{1,4}$
as a double cover of the elliptic curve is given by
\begin{align*}
  y^{2} &= x^{3} - 12x\\
  z^{2} &= (-3x^{4} + 21x^{3} + 192x^{2} + 324x) + (-6x^{2} - 72x - 144)y.
\end{align*}

One can attempt to use \'etale descent. Given a curve $X$ with genus $g > 1$ and plane model $f(x,t) = 0$, an \'etale double cover $Y$ is a curve with an unramified degree $2$ map $\phi : Y \to X$. Since the function field of $Y$, $\Q(Y)$, is a degree two extension of $\Q(X)$, a model for $Y$ can be obtained writing $\Q(Y) = \Q(X)[\sqrt{g(x,t)}]$ for some function $g(x,t) \in \Q(X)$. The model then has the
form
\[
  f(x,t) = 0 \quad y^{2} = g(x,t).
\]
Let $Y_{d}$ be the curve with the equations $f(x,t) = 0$, $y^{2} = dg(x,t)$ and $\pi_{d} : Y_{d} \to X$ the natural map. It is known (see Proposition 5.3.2 of \cite{Skorobogatov} for example) that if $S$ is the set of squarefree integers whose prime factors divide $2$ or the primes of bad reduction of $X$, then $X(\Q) = \bigcup_{d \in S} \pi_{d}(Y_{d}(\Q))$. This reduces the problem of finding the rational points on $X$ to handling the rational points on each $Y_{d}$. Although the curves $Y_{d}$ have genus $2g-1$, they may admit more automorphisms, have maps to lower genus curves, or their Jacobians may satisfy the Chabauty condition. The curve $X$ will admit \'etale double covers defined over $\Q$ if and only if ${\rm Jac}~X$ has a rational $2$-torsion point. 

In the case of the curve $C_{1,4}$, by using the rational $2$-torsion point on $y^{2} = x^{3} - 12x$ to construct a family of \'etale covers. These \'etale covers have genus $7$ and automorphism group $\Z/2\Z \times \Z/2\Z$. There are additional genus $4$ curves that are quotients of these genus $7$ curves. However, we have not found maps to lower genus curves and are unsure how to proceed.

The curve $C_{3,2}$ has genus $2$ and 12 rational points corresponding to $k = 0$, $k = \infty$, $k = 1$, and $k = 135$. A model
for the curve is $y^{2} = x^{6} + 6x^{5} + 39x^{4} + 52x^{3} + 39x^{2}
+ 6x + 1$, and as mentioned in the previous section,
the rational points were provably found by Bianchi and Padurariu in \cite{BP} using quadratic Chabauty.

The curve $C_{3,3}$ has genus $9$ and the quotient by the involution has genus $5$. The Jacobian of this genus $5$ curve appears to be absolutely irreducible, and without maps to lower genus curves it is unclear how to proceed.

We now turn to a discussion of the $D_{1,n}$ and $D_{3,n}$ curves. If $n = 2$,
both of these curves are isomorphic to $\P^{1}$, and $E_{1,k}$ and
$E_{3,k}$ both have a rational point of order $2$ if and only if $k = -\frac{27}{64} a^{3} + (81/64)a - 27/32$ for some rational number $a$.

The curves $D_{1,3}$ and $D_{3,3}$ are both isomorphic to $\P^{1}$ also. We have that $E_{1,k}$ has a rational point of order $3$ if and only if there is a rational number $a$ so that
\[
  k = \frac{-27 (a-1)^{3} (a+1)^{3} (a^{2}+3)}{256a^{2}}
\]
and $E_{3,k}$ has a rational point of order $3$ if and only if there is a rational number $a$ so that
\[
  k = \frac{(a-3)(a+3) (a^{2}+3)^{3}}{256a^{2}}.
\]

The curve $D_{1,4}$ is the elliptic curve $y^{2} = x^{3} - 12x$ with rank $1$ as shown in Section~\ref{proofsec}, Case 4. We have that $E_{1,k}$ has a rational
point of order $4$ on it if and only if there is a rational point $(x,y)$ on $y^{2} = x^{3} - 12x$ for which
\begin{align*}
  k &= \frac{f(x) y + g(x)}{h(x)}, \text{ where }\\
  f(x) &= 5832x^{8} - 151632x^{7} + 1982880x^{6} - 11897280x^{5} + 6438528x^{4}\\
  &+ 330884352x^{3} - 1763596800x^{2} + 2841910272x - 725594112\\
  g(x) &= -729x^{10} + 11664x^{9} - 186624x^{8} + 2869344x^{7} - 33195744x^{6} + 192316032x^{5}\\
  &- 162922752x^{4} - 2831832576x^{3} + 9216052992x^{2} - 6651279360x + 403107840\\
  h(x) &= (x^{2} - 12x - 12)^{6}.
\end{align*}
The curve $D_{3,4}$ is the elliptic curve $y^{2} = x^{3} - 108x$ which has rank zero. The two rational points on $D_{3,4}$ correspond to $k = 0$ and $k = -27/16$. Since the case of $D_{3,4}$ was handled, it is not necessary to consider the cases $D_{3,8}$ and $D_{3,12}$.

The curve $D_{1,5}$ has genus $3$ and plane model $y^{4} = t^{4} -
12t^{3} + 14t^{2} + 12t + 1$. This curve has at least 12 rational
points on it, of which two map to $k = \infty$, two map to $k = 0$,
four map to $k = 864$ and four map to $k = 297/256$. The curve
$D_{1,5}$ has 4 automorphisms over $\Q$ and 16 automorphisms over $K =
\Q(\sqrt{10 + 2 \sqrt{5}})$. There is a degree $2$ map to an elliptic
curve $E$ of rank $2$ defined over $K$, and we attempted to use
elliptic curve Chabauty find the rational points on
$D_{1,5}$. 

Elliptic curve Chabauty is a method developed in \cite{Bruin} to use Chabauty's method on a curve $C$, provided that there is a map $\phi : C \to E$ defined over a number field $K$. The input is a map $\psi : E \to \P^{1}$ that is defined over $K$ with the property that $\phi \circ \psi : C \to \P^{1}$ is defined over $\Q$.
The method essentially applies Chabauty's method to the image of ${\rm Jac}~C$ inside ${\rm Res}_{K/\Q}(E)$. This latter object is a (geometrically disconnected) abelian variety of dimension $[K : \Q]$ whose rank is equal to the rank of $E(K)$. To have any hope for the method to succeed, we must have $\rank E(K) < [K : \Q]$, but this alone is not sufficient, because the image of ${\rm Jac}~C$ inside ${\rm Res}_{K/\Q}(E)$ can be smaller yet.

In case case of $D_{1,5}$, we attempted this. However, it appears that the image of ${\rm Jac} D_{1,5}
\to {\rm Res}_{K/\Q}(E)$ has dimension $2$ and so the Chabauty
condition is not satisfied. If this is the case, the Jacobian of
$D_{1,5}$ has rank $3$ over $\Q$. The curve $D_{1,5}$ is a
double-cover of $y^{2} = t^{4} - 12t^{3} + 14t^{2} + 12t + 1$, which
is isomorphic to the elliptic curve $y^{2} = x^{3} - x^{2} - 13x + 22$
with Mordell-Weil group $\Z/2\Z \times \Z$. Using this, we can create
a family of \'etale double-covers of $D_{1,5}$ using the $2$-isogeny
on the elliptic curve. This gives rise to a family of genus $5$
curves. One such curve has four automorphisms and maps to the
hyperelliptic curve $y^{2} = x^{8} - 8x^{7} + 80x^{6} + 184x^{5} +
430x^{4} + 488x^{3} + 416x^{2} + 200x + 37$. Although this curve has
$16$ automorphisms and many maps to elliptic curves defined over
number fields, the Jacobian of this curve has rank $3$ and so elliptic
curve Chabauty will not work. We have not been able to provably find the
rational points on $D_{1,5}$.

The curve $D_{3,5}$ has the plane model $9y^{4} = t^{4} - 12t^{3} + 14t^{2} + 12t + 1$, and this curve has no points over $\Q_{3}$. Hence, there are no values of $k$ for which $E_{3,k}$ has a point of order $5$. This also implies that it is not necessary to handle $D_{3,10}$.

The curve $D_{1,6}$ has the plane model $y^{4} = t^{4} - 24t^{2} - 48t$
and $D_{3,6}$ has the model $9y^{4} = t^{4} - 24t^{2} - 48t$. In both
cases, there are at least 5 rational points which correspond to $k = 0$ and
$k = \infty$. The 5 rational points generate a subgroup of the Jacobian of
rank 2. These curves have $16$ automorphisms over $\overline{\Q}$,
and there is a map to an elliptic curve with $j = 1728$ of rank $2$
over the number field $L$ defined by $x^{6} + 3x^{4} + 3x^{2} - 3$. (It seems that the image of ${\rm Jac}~D_{1,6}$ inside ${\rm Res}_{L/\Q} E$ is $2$-dimensional and so elliptic curve Chabauty cannot be applied.) The curve $D_{3,6}$
maps to the same elliptic curve. It is worth noting that $D_{1,6}$ and $D_{3,6}$
both have no rational torsion in their Jacobians. In particular, there are
no \'etale double covers defined over $\Q$.

The curve $D_{1,7}$ has the plane model $y^{4} = t^{8} - 4t^{7} - 14t^{6} + 35t^{4} + 56t^{3} + 42t^{2} + 12t + 1$ and the rational points on it were determined in Section~\ref{proofsec}, Case 6. Also in Section~\ref{proofsec}, Case 6, it was shown that $D_{3,7}$ has no local points. 

The curves $D_{1,8}$ has plane model
$y^{4} = t^{8} - 16t^{6} - 32t^{5} + 64t^{3} + 96t^{2} + 64t + 16$. This curve maps to $y^{2} = t^{8} - 16t^{6} - 32t^{5} + 64t^{3} + 96t^{2} + 64t + 16$. This genus
$3$ curve has four automorphisms over $\Q$ and one quotient is the curve $X : y^{2} = -2x^{5} + x^{4} + 8x^{3} - 4x^{2} - 2x + 1$. The Jacobian of this curve has Mordell-Weil group $\Z/2\Z \times \Z$ and Chabauty's method determines that
$X(\Q)$ consists of four points. Pulling these points back to the genus $3$
hyperelliptic curve shows that it has precisely $8$ rational points. Of these
points, the ones that lift to $D_{1,8}$ correspond to $k = 0$ or $k = \infty$. 

The curve $D_{1,9}$ has the model $y^{4} = f(t)$, where
\[
f(t) = t^{12} - 12t^{10} - 28t^{9} - 18t^{8} + 36t^{7} + 114t^{6} + 180t^{5} + 189t^{4} + 128t^{3} + 54t^{2} + 12t + 1.
\]
This polynomial $f(t)$ satisfies $t^{12} f(-1-1/t) = f(t)$ and this gives
rise to an automorphism of $D_{1,9}$ of order $3$. The quotient by this automorphism is the genus $3$ curve $X : -3a^{3}b - 3a^{3} c + b^{3}c + bc^{3} = 0$. The
curve $X$ has a visible automorphism interchanging $b$ and $c$ and the quotient
by this automorphism is the rank one elliptic curve $y^{2} = x^{3} + 9$.
However, $X$ has $48$ automorphisms over $\overline{\mathbb{Q}}$, and there is an
order $2$ automorphism defined over $K = \Q(\sqrt[6]{3})$ for which the
quotient curve is the elliptic curve $y^{2} = x^{3} + (3 - 2 \sqrt{3}) x$.
This elliptic curve has rank $2$ over $K$. We attempted to perform elliptic
curve Chabauty using this setup, but it was not successful. (This suggests
to us that the rank of the Jacobian of $X$ is $3$, although the differences
of rational points on $X$ are torsion in the Jacobian.) We have not been
able to provably determine the rational points on $D_{1,9}$. The curve $D_{3,9}$ has model $9y^{4} = f(t)$, which has no points over $\mathbb{Q}_{2}$.

We did not consider the curves $D_{1,10}$ or $D_{1,12}$. If we were able to provably find the rational points on $D_{1,5}$ and $D_{1,6}$ considering these would not be necessary. 

In summary, there are infinitely many rational $k$ for which there is a rational point on $C_{k}$ whose image on $E_{1,k}$ has order $2$, while there are no $k \ne 0$ for which there is a rational point on $C_{k}$ whose image on $E_{1,k}$ has order $3$. We were not able to handle the
case of points of order $4$ (which does occur for $k = 135$).
The only $k$ for which there is a rational point on $C_{k}$ whose image on $E_{3,k}$ has order $2$ are $k \in \{ 0, \infty, 1, 135 \}$. We were not able to answer the same question for order $3$.

There are infinitely many rational $k$ for which $E_{1,k}$ has a rational point of order $4$. There are finitely many such $k$ for order $5$, $6$, $7$, $8$, $9$, $10$ and $12$, but we are only able to handle orders $7$ and $8$.

There are infinitely many rational $k$ for which $E_{3,k}$ has a rational point of order $3$, while the only $k$ for which $E_{3,k}$ has a rational point of order $4$ are $k \in \{ 0, -27/16 \}$. We are able to show
that $E_{3,k}$ never has a rational point of order $5$ or $7$ or $9$, but we are not able to show the same for order $6$. 

\bibliographystyle{plain}
\bibliography{refs}

\begin{thebibliography}{10}

\bibitem{QC3}
Jennifer Balakrishnan, Netan Dogra, J.~Steffen M\"{u}ller, Jan Tuitman, and Jan
  Vonk.
\newblock Explicit {C}habauty-{K}im for the split {C}artan modular curve of
  level 13.
\newblock {\em Ann. of Math. (2)}, 189(3):885--944, 2019.

\bibitem{QC1}
Jennifer~S. Balakrishnan and Netan Dogra.
\newblock Quadratic {C}habauty and rational points, {I}: {$p$}-adic heights.
\newblock {\em Duke Math. J.}, 167(11):1981--2038, 2018.
\newblock With an appendix by J. Steffen M\"{u}ller.

\bibitem{QC2}
Jennifer~S. Balakrishnan and Netan Dogra.
\newblock Quadratic {C}habauty and rational points {II}: {G}eneralised height
  functions on {S}elmer varieties.
\newblock {\em Int. Math. Res. Not. IMRN}, (15):11923--12008, 2021.

\bibitem{BP}
Francesca Bianchi and Oana Padurariu.
\newblock Rational points on rank 2 genus 2 bielliptic curves in the {LMFDB}.
\newblock Preprint. \url{https://arxiv.org/abs/2212.11635}.

\bibitem{BH}
Yuri Bilu and Guillaume Hanrot.
\newblock Solving {T}hue equations of high degree.
\newblock {\em J. Number Theory}, 60(2):373--392, 1996.

\bibitem{Magma}
Wieb Bosma, John Cannon, and Catherine Playoust.
\newblock The {M}agma algebra system. {I}. {T}he user language.
\newblock {\em J. Symbolic Comput.}, 24(3-4):235--265, 1997.
\newblock Computational algebra and number theory (London, 1993).

\bibitem{Bruin}
Nils Bruin.
\newblock Chabauty methods using elliptic curves.
\newblock {\em J. Reine Angew. Math.}, 562:27--49, 2003.

\bibitem{Kubert}
Daniel~Sion Kubert.
\newblock Universal bounds on the torsion of elliptic curves.
\newblock {\em Proc. London Math. Soc. (3)}, 33(2):193--237, 1976.

\bibitem{Marshall}
David~C. Marshall, Edward Odell, and Michael Starbird.
\newblock {\em Number theory through inquiry}.
\newblock MAA Textbooks. Mathematical Association of America, Washington, DC,
  2007.

\bibitem{Mazur}
B.~Mazur.
\newblock Modular curves and the {E}isenstein ideal.
\newblock {\em Inst. Hautes \'{E}tudes Sci. Publ. Math.}, (47):33--186 (1978),
  1977.
\newblock With an appendix by Mazur and M. Rapoport.

\bibitem{Milne}
J.~S. Milne.
\newblock Abelian varieties.
\newblock In {\em Arithmetic geometry ({S}torrs, {C}onn., 1984)}, pages
  103--150. Springer, New York, 1986.

\bibitem{Silverman}
Joseph~H. Silverman.
\newblock {\em The arithmetic of elliptic curves}, volume 106 of {\em Graduate
  Texts in Mathematics}.
\newblock Springer-Verlag, New York, 1992.
\newblock Corrected reprint of the 1986 original.

\bibitem{SilvermanTate}
Joseph~H. Silverman and John Tate.
\newblock {\em Rational points on elliptic curves}.
\newblock Undergraduate Texts in Mathematics. Springer-Verlag, New York, 1992.

\bibitem{Skorobogatov}
Alexei Skorobogatov.
\newblock {\em Torsors and rational points}, volume 144 of {\em Cambridge
  Tracts in Mathematics}.
\newblock Cambridge University Press, Cambridge, 2001.

\bibitem{Tuitman1}
Jan Tuitman.
\newblock Counting points on curves using a map to {$\bold{P}^1$}.
\newblock {\em Math. Comp.}, 85(298):961--981, 2016.

\bibitem{Tuitman2}
Jan Tuitman.
\newblock Counting points on curves using a map to {$\bold{P}^1$}, {II}.
\newblock {\em Finite Fields Appl.}, 45:301--322, 2017.

\end{thebibliography}

\end{document}